% From uuwendy@senri-i.or.jp Sun Dec 26 21:38:51 1993
% Received: from relay2.UU.NET by sunset.ma.huji.ac.il with SMTP id AA26158
%   (5.65c/HUJI 4.152 for <shlhetal@math.huji.ac.il>); Sun, 26 Dec 1993 21:38:26 +0200
% Received: from spool.uu.net (via LOCALHOST) by relay2.UU.NET with SMTP 
% 	(5.61/UUNET-internet-primary) id AA01708; Sun, 26 Dec 93 14:36:06 -0500
% Received: from iij.UUCP by uucp5.uu.net with UUCP/RMAIL
% 	(queueing-rmail) id 041158.13021; Fri, 24 Dec 1993 04:11:58 EST
% Received: from wincgw1.winc.ad.jp by uucp0.iij.ad.jp (8.6.4+2.2W/2.7W)
% 	id RAA02636; Fri, 24 Dec 1993 17:27:28 +0900
% Received: by wincgw1.winc.ad.jp (5.67+1.6W/4.18:3.7:winc:931216)
% 	id AA05790; Fri, 24 Dec 93 17:27:34 JST
% Received: by synapse.senri-i.or.jp (5.67+1.6W/4.18:3.7:synapse:931220)
% 	id AA23909; Fri, 24 Dec 93 17:30:14 JST
% Received: by japonext.jams.or.jp (NX5.67c/6.4J.6)
% 	id AA06408; Fri, 24 Dec 93 12:48:36 +0900
% Date: Fri, 24 Dec 93 12:48:36 +0900
% From: iseki@jams.or.jp (KIYOSHI ISEKI)
% Message-Id: <9312240348.AA06408@japonext.jams.or.jp>
% Received: by NeXT.Mailer (1.87.2)
% Received: by NeXT Mailer (1.87.2)
% To: shlhetal@math.huji.ac.il
% Subject: no3473-1
% Status: OR
% 

%%%% AMS-TeX file
%%%% with the paper "The number of independent..." by S.Shelah
%%%% submitted to Math. Japonica
%%%% and registered there as No 3473 -I
%%%%%%%%%%%%%   THE FILE STARTS HERE  %%%%%%%%%%%%%%%%%%%%%%%%%
\input amstex
\documentstyle{amsppt}

\def\demo#1{\par \smallskip{\bf #1:}}

\catcode`\@=11

\def\initheadtext#1{\def\initrhead{#1}}
\def\initheadline{%
	\initrhead\hss\headlinefont@\folio}
\headline={\def\chapter#1{\chapterno@. }%
 \def\\{unskip\space\ignorespaces}\headlinefont@
 \iffirstpage@ \global\firstpage@false \initheadline
 \else \ifodd\pageno \rightheadline \else \leftheadline\fi \fi}

\def\output@{\shipout\vbox{%
 \ifrunheads@
	\makeheadline \pagebody
	\else \pagebody \makefootline \fi
 }%
 \advancepageno \ifnum\outputpenalty>\@MM\else\dosupereject\fi}
\catcode`\@=\active
\pagewidth{33truepc}
\pageheight{50truepc}

%%%=================For headline=======================
%\initpageno{000}
\initheadtext{\eightit Math.~Japonica \/ \eightbf 39, \eightrm
		No.~1(1994), 1--5}
\leftheadtext{S. SHELAH}
\rightheadtext{PRODUCT OF BA'S: NU. OF INDEPENDENT}
%===================== Headline end =====================

\document
\topmatter
 \title
THE NUMBER OF INDEPENDENT ELEMENTS
IN THE PRODUCT OF INTERVAL BOOLEAN ALGEBRAS
 \endtitle
 \author Saharon Shelah
 \endauthor
\date Received August 19, 1993 \enddate

 \thanks Partially supported by the Deutsche Forschungsgemeinschaft,
grant Ko 490/7-1. Publication no. 503
 \endthanks
 \address Institute of Mathematics, The Hebrew University of
Jerusalem, Israel; Department of Mathematics, Rutgers University, New
Brunswick, NJ USA
 \endaddress
 \abstract We prove that in the product of $\kappa$ many Boolean
algebras we cannot find an independent set of more than $2^\kappa$
elements solving a problem of Monk (earlier it was known that we
cannot find more than $2^{2^\kappa}$ but can find $2^\kappa$).
 \endabstract
\endtopmatter

\define\pbf{\par\bigpagebreak\flushpar}
\define\un{\underbar}

\define\buildef{\buildrel \text{def} \over =}
\define\uhr{\upharpoonright}

\underbar{\S 0 Introduction}

In his systematic investigation of cardinal invariants of Boolean  
Algebras,
Monk ([M], problem 26, p. 71 or p. 146)
 has raised the following question, where we define:

\noindent
\demo{Definition}  For a Boolean algebra $B$, $\text{Ind}(B)$ is the  
supremum
of $|X|$, for $X \subseteq B$ independent, which means that for  
distinct
$a_1, \dots, a_n \in X$ and non trivial Boolean term $\tau(x_1,  
\cdots, x_n)
$ we have $ B \models\tau(a_1 \cdots, a_n  ) \neq 0$.

Interval Boolean algebras are defined in 1.2 below (non trivial is  
defined
in 1.5 below).
\pbf
\demo{Problem}
If $A_i$ is a non-trivial interval algebra for each $i \in I$, where  
$I$ is
infinite, is $\text{Ind}(\Pi_{i \in I} A_i) = 2^{|I|}$?   
Equivalently, is it
true that for every infinite cardinal $\kappa$ there is no linear  
order $L$
and sequence $\langle x_\alpha$ : $ \alpha < (2^\kappa)^+ \rangle $  
with the
following properties?
\item{(1)} For all $\alpha < (2^\kappa)^+$, $x_\alpha$ is a sequence
$\langle x_{\alpha, \beta}: \beta < \kappa \rangle$ such that for  
$\beta <
\kappa$, $x_{\alpha, \beta}$ is a finite collection of half-open  
intervals
of $L$.
\item{(2)} For all finite disjoint $\Gamma, \nabla \subseteq
(2^\kappa)^+ $ there is a $\beta < \kappa$ such that
$$\mathop{\bigcap}\limits_{\alpha \in \Gamma} (\bigcup x_{\alpha,  
\beta})
\cap \mathop{\bigcap}\limits_{\alpha \in \nabla} (L \setminus
\bigcup x_{\alpha, \beta}) \neq 0.$$

It was known that $2^{2^{\kappa}} \ge \text{Ind}(\underset \zeta <  
\kappa
\to{\Pi} B_\zeta) \ge 2^\kappa$, and it was felt that the answer to  
the
question is independent of ZFC.  Monk phrased some weaker related  
questions
of interest to set theorists on this see Shelah and Soukup [ShSo376].

But not all problems in set theory are independent of set theory and  
as we
can see below, the main point is to get those with an answer and here  
we
even get a reasonable one, so we discuss to some length how simple.   
Those
results are essentially best possible as will be shown elsewhere.  In  
his
lecture in Jerusalem, Monk raised the question again, for which we  
thank
him, and Mati Rubin had enough interest in the solution to rewrite it
beautifully for which I thank him and the reader should too.  Lately  
I have
learned that: Just and Weese had gotten a restricted positive answer
(assuming $B_i = B(I_i), \; I_i$ a well order), and Rubin [R, lemma  
5.4]
proved: if $\lambda$ is a regular cardinal and $\{a_\alpha | \alpha <
\lambda\}$ is a sequence of elements of an interval algebra $B$, then
$\{a_\alpha |\alpha < \lambda \} $ contains a semihomogeneous  
subsequence of
length $\lambda$.

\pbf
\un{Notation}: The operations in a Boolean algebra are denoted by  
$\cdot$
(product, intersection) $+$ (addition, union), $-$ (complement), and
$\triangle$ is the symmetric difference.
\pbf
\underbar{\S 1 The main result}

\proclaim{Theorem 1.1} Let $\kappa$ be an infinite cardinal, and for  
every
$\zeta < \kappa$ let $B_\zeta$ be an interval Boolean algebra.  Then
$\text{Ind}(\underset \zeta < \kappa\to{\Pi} B_\zeta) = 2^\kappa$.

Moreover, there is $n \in \omega$ and a nontrivial Boolean term  
$\tau(x_{0,
\cdots,} x_{n-1})$ such that for every $\{ a_\alpha  
|\alpha<(2^\kappa)^+
\} \subseteq \underset \zeta< \kappa \to{\Pi} B_\zeta$,  there are
$ \alpha_0<\cdots< \alpha_{n-1} $ such that $\tau(a_{\alpha_0,  
\cdots,}
a_{\alpha_{n-1}}) = 0$. In fact $\tau$ can be taken to be
$\tau(x_0, \cdots, x_5) \buildrel \text{def} \over = x_0 \cdot x_1  
\cdot
(-x_2) \cdot (-x_3) \cdot x_4 \cdot (-x_5)$.
\endproclaim
We need some notations and definitions.

\demo{Definition 1.2} If $\langle I, < \rangle$ is a linear ordering  
let
$I^+ \buildef I \cup \{ - \infty, \infty \}$.  We assume that  
$-\infty,
\infty \notin I$, and define the linear ordering $<^+$ on $I^+$ in  
the
obvious way.  For $s, t \in I^+$, let $(s, t) \buildef \{ x \in I| s  
<^+ x
<^+ t \}$, $[s, t) \buildef \{ x \in I | s \le^+ x <^+ t \}$ etc. Let  
$B(
\langle I, < \rangle)$ be the subalgebra of $\Cal P(I)$ generated by  
$\{[s,
t) | s, t \in I^+ \}$.  We abbreviate $<^+$ by $<$ and $B(\langle I,  
<
\rangle )$ by $B(I)$.  A Boolean algebra of this form is called an  
interval
Boolean algebra.  Let $I^* = \{ x \in I |$ there is $y \in I$ such  
that $ y
< x \} \cup \{ - \infty, \infty \}$.  Every $a \in B(I)$ has a unique
representation of the form $a =\underset i <n \to{\cup} [s_{2i},  
s_{2i+1}
)$, where $ n \ge 0$, $s_{0, \cdots,} s_{2n-1} \in I^*$ and $s_0 <  
s_1 <
\ldots < s_{2n-1}$. We denote $\sigma_a^- = \{ s_{0, \cdots,}  
s_{2n-1} \} $,
$\sigma_a = \sigma_a^- \cup \{ - \infty, \infty \}, n_a = |\sigma_a  
|$ and
$\vec \sigma_a = \langle s'_{0, \cdots,}  s'_{n_a - 1} \rangle$,  
where $\{
s'_i | i < n_a \} = \sigma_a$ and $s'_0 < \dots < s'_{n_a -1}$.

\noindent
\demo{Definition 1.3} Let $\langle I, < \rangle$ be a linear  
ordering,
$S$ be a set of ordinals, and $\vec a= \{ a_\alpha|
\alpha \in S \} \subseteq B(I)$.

$(a) \;  \vec a$ is homogeneous, if:

\noindent
\ \ \ \ \ \ \ \ (1) there is $k \in \omega$ such that for every  
$\alpha \in
S, \; \; |\sigma_{a_\alpha}| = k$;

\noindent
\ \ \ \ \ \ \ \ (2) for every $\alpha, \beta \in S$,  
$\sigma^-_{a_\alpha}
\cap \{ - \infty, \infty \} = \sigma^-_{a_\beta}\cap \{ -  \infty,  
\infty
\}$, and

\noindent
\ \ \ \ \ \ \ \ (3) for every $\alpha < \beta $ in $S$ there is  
$\ell_{\alpha,
\beta} < n_{a_\alpha} -1$ such that for every

\noindent
\ \ \ \ \ \ \ \ \ \ \  $s \in \sigma_{a_\beta}  
\setminus\{-\infty,\infty\}$,
$\vec \sigma_{a_\alpha} (\ell_{\alpha, \beta}) < s <  \vec
\sigma_{a_{\alpha}} (\ell_{\alpha, \beta} + 1)$.

(b) $\vec a$ is semi-homogeneous, if there are $-\infty = t_0 <
\dots < t_k = \infty$ in $I^+$ such that for every $\ell < k
$ we have: $\{ a_\alpha\cap [t_\ell, t_{\ell + 1}) | \alpha\in S \} $  
is
homogeneous in $B(I) \uhr [t_\ell, t_{\ell + 1})$.  (Note that $B(I)  
\uhr
[t_\ell, t_{\ell + 1})$ is the interval algebra of the linear  
subordering of
$I$ whose universe is $[t_\ell, t_{\ell + 1}) \; )$. Now $ \{ t_{0,  
\cdots,
} t_k \}$ is called a partitioning set for $\vec a$.

If $\{ B_\zeta|\zeta <\kappa\}$ is a family of $BA's$, then the  
members of
$\underset \zeta<\kappa\to{\Pi} B_\zeta$ are denoted by $\langle  
a_\zeta|
\zeta< \kappa\rangle$.  So for every $\zeta < \kappa$, $a_\zeta \in
B_\zeta$.

The following claim is our main lemma.
\proclaim{Lemma 1.4} Let $\kappa$ be an infinite cardinal.
For every $\zeta <\kappa$ let $B_\zeta$ be an interval Boolean  
algebra.  Let
$\vec a = \left\{ a_\alpha | \alpha <(2^\kappa)^+\right\}
\subseteq \underset \zeta<\kappa\to{\Pi} B_\zeta$, and denote  
$a_\alpha =
\langle a_{\alpha,_\zeta}| \zeta < \kappa\rangle$. Then there is $S
\subseteq (2^\kappa)^+$ of cardinality $(2^\kappa)^+$ such that for  
every
$\zeta <\kappa $ we have: $ \{ a_{a,\zeta} | \alpha \in S \}$ is
semi-homogeneous.  In fact, $S$ can be taken to be a stationary  
subset of
$(2^\kappa)^+$.
\endproclaim
\demo{Proof} Let $\langle I, < \rangle$ be a linear ordering, $a \in  
B(I)$
and $A \subseteq I$.  A set $C \subseteq A \cup \{ - \infty, \infty  
\}$ is
called an $A$-partition of $a$, if: $\{ - \infty, \infty \} \subseteq  
C$,
$\sigma_a \cap A\subseteq C$, and for every $\ell <n_a -1$: if $(\vec
\sigma_a (\ell), \vec \sigma_a (\ell+1))
\cap A \neq \emptyset$,
then $(\vec \sigma_a (\ell), \vec \sigma_a (\ell + 1)) \cap C \neq
\emptyset$.

Let $ \{ \langle I_\zeta, <_\zeta \rangle | \zeta < \kappa\}$ be a  
family of
linear orderings, such that $B_\zeta = B(I_\zeta)$, let $\lambda =
(2^\kappa)^+$ and $\{ a_\alpha | \alpha < \lambda \} \subseteq
{\underset{\zeta < \kappa}\to{\Pi}} B_\zeta$.  We may assume that  
$|I_\zeta|
= \lambda$ for every $\zeta < \kappa$, and hence we may further  
assume that
$I_\zeta = \lambda$.  For $s, t \in \lambda$ and $\zeta < \kappa$, we  
use
$[s, t)_\zeta, (s, t)_\zeta$ etc. to denote intervals of $ <_\zeta$  
whose
endpoints are $s$ and $t$.  For every $\alpha < \lambda$ let  
$a_\alpha =
\langle a_{\alpha, \zeta} |\zeta < \kappa \rangle$.  For every  
$\alpha <
\lambda$, $\zeta < \kappa$ and $\ell \le n_{a_{\alpha, \zeta}} - 1$,  
let
$\sigma_{\alpha, \zeta} = \sigma_{a_{\alpha, \zeta}},
\vec \sigma_{\alpha, \zeta} = \vec \sigma_{a_{\alpha,\zeta}}$ and
$s_{\alpha, \zeta, \ell} = \vec \sigma_{a_{\alpha, \zeta}} (\ell)$.

We next perform on the sequence $\{ a_\alpha | \alpha < \lambda \}$  
several
steps of cleaning.  Let $S_0 = \{ \alpha < \lambda |  
\text{cf}(\alpha) =
\kappa^+ \}$.

\underbar{Step 1}: By partitioning $S_0$ into $4^\kappa < \lambda$   
sets, we
obtain a stationary set $S_1 \subseteq S_0$ such that for every  
$\alpha,
\beta \in S_1$ and $\zeta <\kappa$: \ \ $\sigma^-_{\alpha, \zeta}  
\cap \{
-\infty, \infty \} = \sigma^-_{\beta, \zeta} \cap \{ -\infty, \infty  
\}$.

\underbar{Step 2}: For every $\alpha < \lambda$ let $\vec n_\alpha =
\langle n_{a_{\alpha, \zeta}} | \zeta < \kappa \rangle $.
For every $\vec n \in \; {}^\kappa \omega$ let $ S^1_{\vec n} = \{  
\alpha
\in S_1 | \vec n_\alpha = \bar n \}$.  So $ \{S^1_{\vec n} | \; \vec  
n \in
\; ^{\kappa}\omega \}$ is a partition of $S_1$ into $\le 2^\kappa <  
\lambda
$ sets.  Hence for some $\vec n \in \; ^{\kappa}\omega, \; S_2  
\buildef
S^1_{\vec n}$ is stationary.  Let $\vec n = \langle n_\zeta | \zeta <
\kappa\rangle$.

\underbar{Step 3}: For every $\alpha \in S_2$ let
$C_\alpha \subseteq \alpha \cup \{ -\infty, \infty \}$ have the  
following
properties:

\noindent
(1) \ $|C_\alpha| \le \kappa$; and

\noindent
(2) for every $\zeta <\kappa, C_\alpha$ is an $\alpha$-partition of
$a_{\alpha, \zeta}$ with respect to $<_\zeta$.

Since $S_2 \subseteq \{ \alpha < \lambda | \text{cf}(\alpha) =  
\kappa^+
\}$, clearly $\sup(C_\alpha \setminus \{ - \infty, \infty  
\})<\alpha$.  By
Fodor's theorem and the fact that $(2^\kappa)^\kappa= 2^\kappa$,  
there is a
stationary set $S_3 \subseteq S_2$, $\delta < \lambda$ and $C  
\subseteq
\delta \cup \{ -\infty, \infty \}$ such that for every $\alpha \in  
S_3,
C_\alpha = C$.

\underbar{Step 4}: By partitioning $S_3$ into $\le 2^\kappa$ sets we  
obtain
a stationary set $S_4 \subseteq S_3$ and a system  
$T_\zeta=\{t_{\zeta, 0},
\dots, t_{\zeta, m_\zeta } \}\subseteq C$ for $ \zeta <\kappa$, such  
that
for every $\alpha, \beta \in S_4$ and $\zeta < \kappa$:

\noindent
(1)\ \ \ $-\infty =t_{\zeta, 0} <_\zeta t_{\zeta, 1} < _\zeta\dots  
<_\zeta
t_{\zeta, m_{\zeta}}=\infty;$

\noindent
(2)\ \ \ $T_\zeta$ is a $(C \setminus \{ -\infty, \infty  
\})$-partition of
$a_{\alpha, \zeta}$;

\noindent
(3)\ \ \ for every $\ell \le n_\zeta -1:$ if $s_{\alpha,  \zeta,  
\ell} \in
T_\zeta$, then  $s_{\beta, \zeta, \ell} = s_{\alpha, \zeta, \ell}$,  
and

\noindent
(4)\ \ \ for every $\ell <n_\zeta -1,\; \; T_\zeta \cap  
\left(s_{\alpha,
\zeta, \ell}, s_{\alpha, \zeta, \ell+1}\right)_\zeta = T_\zeta \cap
\left(s_{\beta, \zeta,\ell,} s_{\beta, \zeta, \ell +1}\right)_\zeta$.

\underbar{Step 5}: Let $F \subseteq \lambda$ be a closed and  
unbounded
set such that for every $\gamma \in F$ and $\alpha < \gamma$ and
$\zeta<\kappa $ we have $ \sigma_{a_{\alpha,\zeta}}
\setminus \{ - \infty, \infty \}
\subseteq \gamma$.  Let $S = S_4 \cap F$.
\pbf
We shall show that for every $\zeta<\kappa$:
\par\noindent
$(\ast)_\zeta$\ \ \ \ \ $\vec a^\zeta \buildef \{ a_{\alpha, \zeta} |  
\alpha
\in S \}$ is semi-homogeneous, and $T_\zeta$ is a partitioning set  
for $\vec
a^\zeta$.
\pbf
Let $m < m_\zeta$, and we show that $\{ a_{\alpha, \zeta} \cap  
[t_{\zeta,
m}, t_{\zeta, m+1})_\zeta|\alpha \in S \}$ is homogeneous in $B_\zeta
\uhr[t_{\zeta, m,} t_{\zeta, m+1})_\zeta$.  So for the rest of the  
proof of
1.4, we fix $\zeta$ and $m$.  Let $a'_\alpha = a_{\alpha, \zeta} \cap
[t_{\zeta, m}, t_{\zeta, m+1})_\zeta$.  For every $\alpha \in S$:

$T_\zeta $ is a $(C \setminus \{ -\infty, \infty \}) $-partition of
$a_{\alpha, \zeta}$ and $C$ is an $\alpha$-partition of $a_{\alpha,
\zeta}$.
\par\noindent
Hence
\par \flushpar
(i) \ \ \ \ \ for every $\alpha \in S$: $T_\zeta$ is an  
$\alpha$-partition
of $a_{\alpha, \zeta}$.

Let $I' = [ t_{\zeta, m}, t_{\zeta, m+1})_\zeta$ and $ <'=<_\zeta
\upharpoonright I'$ and $B'=B(I')$. For every $\alpha \in S$ let  
$n'_\alpha
= n_{a'_\alpha}, \sigma'_\alpha = \sigma_{a'_\alpha} $ and $\langle
s'_{\alpha, 0, \cdots, } s'_{\alpha,n'_\alpha - 1} \rangle = \vec
\sigma_{a'_{\alpha}},$ where the representation is taken with respect  
to
$B'$.  For $\alpha \in S: $ if $t_\zeta, {}_m \neq - \infty $ let
$\ell^\alpha_0$ be such that $t_{\zeta, m} \in [s_{\alpha, \zeta,
\ell^\alpha_0}, s_{\alpha, \zeta, \ell^\alpha_0+1})_\zeta$ and
$\ell^\alpha_0 = - \infty$ otherwise,  and if $t_{\zeta, m+1} \neq  
\infty$,
let $\ell^\alpha_1$ be such that $t_{\zeta, m+1} \in [s_{\alpha,  
\zeta,
\ell^\alpha_1}, s_{\alpha,\zeta, \ell^\alpha_1 +1})_\zeta $ and $
\ell^\alpha_1 = \infty$ otherwise. By conditions (3) and (4) in step  
4, for
every $\alpha, \beta \in S: \; \; \ell^\alpha_0 = \ell^\beta_0$, and
$\ell^\alpha_1 = \ell^\beta_1$. It follows that for every $\alpha,  
\beta \in
S: n'_\alpha = n'_\beta$  and $\{ - \infty, \infty \} \cap
\sigma^-_{a'_{\alpha}} = \{ -\infty, \infty \} \cap  
\sigma^-_{a'_{\beta}}$.
Let $n' = n'_\alpha$.  This means that $\{ a'_\alpha | \alpha \in S  
\}$
satisfies conditions (1) and (2) in 1.3(a).

Let $\alpha < \beta$ be in $S.$  By the choice of $F$ and $S$, (ii)  
holds:
\par\noindent
(ii)\ \ \ \ \ $\sigma'_\alpha \setminus \{ - \infty, \infty\}  
\subseteq
\sigma_{a_{\alpha, \zeta}} \setminus \{ -\infty, \infty \} \subseteq  
\beta$.

By (i), $T_\zeta $ is a $\beta$-partition of $a_{\beta, \zeta},$ and  
hence
$$
\sigma'_\beta \cap \beta\setminus\{ - \infty, \infty \}\subseteq
\sigma_{a_{\beta, \zeta}} \cap \beta \cap (t_{\zeta, m}, t_{\zeta,
m+1})_\zeta = \emptyset.
$$
It follows that for every distinct $\alpha, \beta
\in S$, $(\sigma'_\alpha \setminus \{ - \infty, \infty \} ) \cap
(\sigma'_\beta \setminus \{ -\infty, \infty \}) = \emptyset$.

Suppose by contradiction that for some $\alpha < \beta \in S, k \le  
n'-1$
and $\ell < n'-1$ we have: $s'_{\beta, \ell}, s'_{\beta,\ell+1}  
\notin
\{ - \infty, \infty \} $ and $s'_{\beta, \ell} <_\zeta s'_{\alpha, k}
<_\zeta s'_{\beta,\ell+1}, \; \; s'_{\alpha, k} \notin \{ -\infty,  
\infty
\}$, and hence by (ii), $s'_{\alpha, k} <\beta$. So $(s'_{\beta,  
\ell,}
s'_{\beta, \ell + 1})_\zeta \cap \beta \neq \emptyset.$
Since $s'_{\beta, \ell}, s'_{\beta_{\ell + 1}} \notin \{ - \infty,  
\infty \}$,
 there is $\ell_1 < n_\zeta -1$, such that $s'_{\beta, \ell} =  
s_{\beta, \zeta,
\ell_1} $ and $s'_{\beta, \ell + 1} = s_{\beta, \zeta, \ell_1 +1}$.
So $(s_{\beta, \zeta, \ell_1}, s_{\beta, \zeta, \ell_1 +1})_\zeta  
\cap \beta
\neq \emptyset$. Since $T_\zeta$ is a $\beta$-partition for  
$a_{\beta,
\zeta}$, we have $T_\zeta \cap (s_{\beta, \zeta, \ell_1,} s_{\beta,  
\zeta,
\ell_1 +1})_\zeta \neq \emptyset$. A contradiction. This shows that
$\{a'_\alpha | \alpha \in S \}$ satisfies condition (3) in 1.3(a).
So $\{ a'_\alpha | \alpha \in S \} $  is homogeneous, and hence
$\{a_{\alpha, \zeta} | \alpha \in S \} $ is semi-homogeneous.
\hfill $\square$
\pbf
\demo{Definition 1.5}
Let $\tau(x_1, \cdots, x_n)$ be a Boolean term.  $\tau$ is  
nontrivial, if
for some Boolean algebra $ B$ and $ a_1, \cdots, a_n \in B$ we have $  
\tau
(a_1, \cdots, a_n) \neq 0.$
\proclaim{Lemma 1.6} There are nontrivial Boolean terms
$\tau_1(x_0, x_1, x_2), \; \tau_2 (x_0, x_1, x_2), \; \tau_3 (x_1,  
x_2)$ and
$\tau_4 (x_1, x_2)$ such that for every interval Boolean algebra $B$  
and a
homogeneous sequence $\{ a_i | i < 3 \} \subseteq B$ we have: $ B  
\models
\overset 4 \to {\underset i = 1 \to {\bigvee}} (\tau_i = 0) [v]$,  
where $v$
is the assignment that takes each $x_i$ to $a_i$.
\endproclaim
\demo{Proof} Let $\tau_1 = x_0 \cdot x_1 \cdot (-x_2), \tau_2 =  
(-x_0) \cdot
(-x_1) \cdot x_2, \tau_3 = x_1 \cdot x_2 $ and $\tau_4 = (-x_1) \cdot
(-x_2)$.  Let $B = B(\langle I, < \rangle )$ be an interval algebra  
and $\{
a_i | i < 3 \} $ be homogeneous.
\par \noindent
Let $\vec\sigma_{a_i} = \langle s^i_0, \cdots, s^i_{n-1}\rangle.$
For every $i < j <3$ let $\ell_{i, j} < n-1$ be such that for
every $s \in \sigma_{a_j} \setminus \{ -\infty, \infty \}, \;
s^i_{\ell_{i, j}} <s<s^i_{\ell_{i, j}+1}$ (it exists - see Def  
1.3(a)(3)).

\noindent
\underbar{Case 1} \   $\ell_{1,2} \neq 0$ and $\ell_{1, 2} + 1 \neq  
n-1$.
 It follows that $\ell_{0, 2} = \ell_{0, 1}$, and so $a_1 \triangle  
a_2
\subseteq [s^0_{\ell_{0, 1}}, s^0_{\ell_{0, 1} +1})
\buildef J$.  Either (i) $J \subseteq a_0$, or (ii) $J\cap a_0 =  
\emptyset.$
Now (i) implies that $B \models (\tau_2 = 0)[v]$, and (ii) implies  
that $B
\models (\tau_1 = 0) [v]$.

\noindent
\underbar{Case 2} \  $\ell_{1, 2} = 0$ or $\ell_{1, 2} +1 = n - 1$.   
W.l.o.g.
$ \ell_{1, 2} = 0$. Let $J_1 = [-\infty, s^1_1)$ and $J_2 =  
[s^2_{n-2},
\infty)$. So:
\par\noindent
(1) \ \ \ \ $ J_1 \cup J_2 = I$;
\par\noindent
(2) \ \ \ \ $ J_1 \subseteq a_1 $ or $ J_1 \subseteq -a_1$ ; and
\par\noindent
(3) \ \ \ \ $J_2\subseteq a_2$ or $J_2 \subseteq -a_2$.
\par\noindent
Hence there are four possibilities:
$a_1 + a_2 = 1, a_1 + (-a_2) = 1, (-a_1) + a_2 = 1$ or $ (-a_1) +  
(-a_2) =
1$.  It is now trivial to check that $B \models \overset 4  
\to{\underset i=1
\to{\bigvee}} (\tau_i = 0) [v]$ in all the above subcases. (Note, if
$\ell_{1, 2} + 1 = n -1$.  Use $J'_1 = [s^1_{n-2}, \infty), J_2' =  
[-\infty,
s^2_1)$).
\hfill $\square $
\demo{1.7 Proof of theorem 1.1} Let
$$ \tau(x_0, \cdots, x_{11}) = \tau_1(x_0, x_1, x_2) \cdot  
\tau_2(x_3, x_4,
x_5)\cdot\tau_3(x_7, x_8) \cdot \tau_4(x_{10}, x_{11}).$$
For every $\zeta <\kappa$ let $B_\zeta=B(\langle I_\zeta, <_\zeta  
\rangle )$
be an interval algebra. Let $\lambda=(2^\kappa)^+,$ and $\{ a_\alpha  
|
\alpha < \lambda \} \subseteq \prod_{\zeta <\kappa} B_\zeta$, where
$a_\alpha = \langle a_{\alpha, \zeta} | \zeta <\kappa\rangle$. By  
lemma 1.4
we may assume that for every $\zeta<\kappa$, $\{ a_{\alpha, \zeta} |
\alpha < \lambda \} $ is semi-homogeneous with the partitioning  
sequence
$ \langle t_{\zeta, 0}, \cdots t_{\zeta, m_\zeta} \rangle$.
For every $\zeta < \kappa$ and $m < m_\zeta$ let $B_{\zeta, m} =  
B_\zeta
\uhr [t_{\zeta, m}, t_{\zeta, m+1})$ and $a_{\alpha, \zeta, m} =  
a_{\alpha,
\zeta} \cap [t_{\zeta, m}, t_{\zeta, m+1})$.
Hence $B = \prod \{ B_{\zeta, m} | \zeta < \kappa, m <m_\zeta \}$ and
$a_\alpha = \langle a_{\alpha, \zeta, m} | \zeta< \kappa, m<m_\zeta
\rangle$, and for every $\zeta <\kappa$ and $m <m_\zeta, \{  
a_{\alpha,
\zeta, m} | \alpha < \lambda \} $ is homogeneous in $B_{\zeta, m}$.
So by renaming $\{ B_{\zeta, m} | \zeta < \kappa, m < m_\zeta \}$ as
$\{ B'_\zeta|\zeta<\kappa\}$ and $\{ a_{\alpha, \zeta, m} | \zeta  
<\kappa, m
<m_\zeta \}$  as $\{ a'_{\alpha, \zeta} | \zeta <  \kappa\}$, we may  
assume
that for every  $\zeta < \kappa, \; \;  \{ a_{\alpha, \zeta} | \zeta<
\kappa\}$ is homogeneous in $B_\zeta$. For every $\alpha < \beta <  
\lambda $
and $\zeta  < \kappa$ let $\vec \sigma_{a_{\alpha,\zeta}} = \langle
s^{\alpha, \zeta}_0, \cdots, s^{\alpha, \zeta}_{n_{\zeta  
-1}}\rangle$, and
let $\ell = \ell^\zeta_{\alpha, \beta} < n_\zeta$ be such that for  
every $s
\in \sigma_{a_{\beta, \zeta}} \setminus \{ -\infty, \infty \}$ we  
have $
s^{\alpha, \zeta}_\ell < s < s^{\alpha, \zeta}_{\ell + 1}$.
Let $\vec \ell_{\alpha, \beta} = \langle \ell^\zeta_{\alpha, \beta} |
\zeta < \kappa \rangle$.

There are four triples $ \alpha_0 < \cdots < \alpha_{11} < \lambda$,
such that for every $i =1, 2, 3 \; \vec \ell_{\alpha_0 \alpha_1} =
\vec \ell_{\alpha_{3i}, \alpha_{3i +1}}$ and $\vec \ell_{\alpha_1,  
\alpha_2}
 = \vec \ell_{\alpha_{3i+1}, \alpha_{3i+2}}$. So for every Boolean  
term
$\tau(x_0, x_1, x_2), \zeta<\kappa$ and $i = 1, 2, 3$
we have: $B_\zeta \models \tau (a_{\alpha_0, \zeta}, a_{\alpha_1  
\zeta},
a_{\alpha_2, \zeta})= 0$ iff $B_\zeta \models \tau(a_{\alpha_{3 i,  
\zeta}},
a_{\alpha_{3i+1, \zeta}}, a_{\alpha_{3i+2}, \zeta}) = 0$. Since  
$\overset 4
\to{\underset i = 1 \to{\bigvee}} (\tau_i (a_{\alpha_0, \zeta,}
a_{\alpha_1, \zeta}, a_{\alpha_2, \zeta}) = 0)$ holds in $B_\zeta$  
clearly,
$$
\tau(a_{\alpha_0, \zeta}, \cdots, a_{\alpha_{11}, \zeta}) =
\prod^4_{i=1} \tau_i(a_{\alpha_{3i, \zeta}}, a_{\alpha_{3i+1,  
\zeta}},
a_{\alpha_{3i+2, \zeta}}) =0.
$$
That is, every coordinate of
$\tau\left(a_{\alpha_0 , \cdots}, a_{\alpha_{11}}\right)$ is equal to
$0$.  So $\tau(a_{\alpha_0, \cdots}, a_{\alpha_{11}}) = 0.$ Note that
$\tau$ really has only eight variables.

\noindent
\demo{1.8 Proof with the shorter term}
In order to show that $\tau$ can be taken to be $x_0 \cdot x_1 \cdot  
(-x_2)
\cdot (-x_3) \cdot x_4 \cdot (-x_5)$, we first notice that for every  
$\alpha
< \beta < \gamma$: if $ \vec \ell_{\alpha, \beta} = \vec  
\ell_{\alpha,
\gamma}$, then for every $\zeta < \kappa$: $(a_{\beta, \zeta}  
\triangle
a_{\gamma, \zeta}) \cdot a_{\alpha, \zeta} = 0$, or $(a_{\beta,  
\zeta}
\triangle a_{\gamma, \zeta})\cdot -a_{\alpha,\zeta} = 0,$
moreover from the value of $\vec \ell_{\alpha, \beta} = \vec  
\ell_{\alpha,
\gamma}$: (but not $\vec \ell_{\beta,\gamma}$) we can compute an  
equation,
which is one of those two and holds (possibly both holds).

For every $\alpha < \lambda$ let $L_\alpha=\{ \vec \ell | \quad
| \{ \beta > \alpha | \vec \ell_{\alpha, \beta } = \vec \ell \} | =
\lambda \}$.
So since the number of possible $\vec \ell_{\alpha, \beta}$'s is  
$2^\kappa$,
$L_\alpha \neq \emptyset$. For every $\vec \ell$ let $\Lambda_{\vec  
\ell}
= \{ \alpha | \vec \ell \in L_\alpha \}$. So for some $\vec \ell^0$,  
$ |
\Lambda_{\vec \ell^0} | = \lambda$. Let $\alpha_0 < \cdots <  
\alpha_5$
be such that $\vec \ell_{\alpha_0, \alpha_1} = \vec \ell_{\alpha_0,  
\alpha_2}
= \vec \ell_{\alpha_3, \alpha_4} = \vec \ell_{\alpha_3, \alpha_5} =  
\vec
\ell^0$. (We can demand that $\vec \ell_{\alpha_1, \alpha_2} = \vec
\ell_{\alpha_4, \alpha_5}$, but it is not needed). Let $\zeta <  
\kappa$.
Then either
$$\overset 1 \to{\underset i=0\to{\bigwedge}}
(a_{\alpha_{3i+1},\zeta} \triangle a_{\alpha_{3i+2}, \zeta}) \cdot
a_{\alpha_{3i}, \zeta} = 0$$
 or
$$\overset 1 \to{\underset i=0 \to{\bigwedge}}
(a_{\alpha_{3i+1}, \zeta} \triangle a_{\alpha_{3i+2}, \zeta}) \cdot  
(-a_{3i,
\zeta}) = 0.$$
It follows that $(a_{\alpha_1, \zeta} \triangle a_{\alpha_2,
\zeta}) \cdot a_{\alpha_0, \zeta} \cdot (a_{\alpha_4,\zeta}\triangle
a_{\alpha_5,\zeta})\cdot(-a_{\alpha_3, \zeta}) = 0$. So  
$\tau(a_{\alpha_0,
\zeta}, \cdots, a_{\alpha_5, \zeta}) = 0.$  Hence $\tau$ is as  
required.
\hfill $\square $

\noindent
\proclaim{Claim 1.9} In 1.1 we can use $(x_0 \triangle x_1) \cdot x_2  
\cdot
(x_3 \triangle x_4) \cdot (-x_5)$
\endproclaim

\noindent
\demo{Proof} As above but in 1.8 demand also
$\vec \ell_{\alpha_1, \alpha_2} = \vec \ell_{\alpha_4, \alpha_5} $.

\noindent
\proclaim{Claim 1.10} Let $\kappa$ be an infinite cardinal, and for  
$\zeta
< \kappa$ let $B_\zeta$ be an interval Boolean algebra.  If $a_\alpha  
\in
\underset \zeta < \kappa \to{\Pi} B_\zeta$ for $\alpha <  
(2^\kappa)^{++}$,
then for some $\alpha_0 < \alpha_1 < \alpha_2 < \alpha_3 <  
(2^\kappa)^{++}$
we have $$(a_{\alpha_0} \triangle a_{\alpha_1}) \cdot (a_{\alpha_2}
\triangle a_{\alpha_3}) = 0.$$
\endproclaim
\noindent
\demo{Proof} Similar only in 1.7, 1.8  we use the easy fact (which  
follows
from Erd\"os--Rado)
\par\noindent
$(*)$\ \ \ \ \ if $c$ is a two place function from $(2^\kappa)^{++}$  
to
$\kappa$
\par\noindent
\ \ \ \ \ \ \ then for some $\alpha_0 < \alpha_1 < \alpha_2 <  
\alpha_3$ we
have $c(\alpha_0, \alpha_2) = c(\alpha_0, \alpha_3) = c(\alpha_1,  
\alpha_2)
= c(\alpha_1, \alpha_3) $.

\enddocument

\bye